\documentclass[reqno]{amsart}
\usepackage{marginnote}
\usepackage[marginparwidth=.7in]{geometry}

\usepackage{amscd,amsthm,amsfonts,amssymb,amsmath}
\usepackage{amsmath,tikz-cd}
\usepackage[colorlinks=true]{hyperref}
\usepackage{fullpage}
\usepackage[all]{xy}
\usepackage{mathtools}
\usepackage{mathrsfs}
\usepackage{dsfont}
\usepackage[small,nohug,heads=vee]{diagrams}

\usepackage{hyperref}
\usepackage{color}
\hypersetup{
    colorlinks=true,
    citecolor=black,
    linkcolor=black,
  urlcolor=blue,
}
 
\newtheorem{problem}{Problem}

     \newcommand{\BF}{{\mathbb {F}}}

     \newcommand{\BN}{{\mathbb {N}}}
     \newcommand{\BP}{{\mathbb {P}}}
    \newcommand{\BQ}{{\mathbb {Q}}}

     \newcommand{\BZ}{{\mathbb {Z}}}

    \newcommand{\CO}{{\mathcal {O}}} 
     \newcommand{\CR}{{\mathcal {R}}}

     \newcommand{\CX}{{\mathcal {X}}}

     \newcommand{\fp}{{\mathfrak{p}}}

    \newcommand{\Ad}{{\mathrm{Ad}}}

    \newcommand{\Aut}{{\mathrm{Aut}}}

    \newcommand{\Cor}{{\mathrm{Cor}}}

    \newcommand{\End}{{\mathrm{End}}} 
    
    \newcommand{\Fr}{{\mathrm{Fr}}}
    \newcommand{\Frob}{{\mathrm{Frob}}}

    \newcommand{\Gal}{{\mathrm{Gal}}} \newcommand{\GL}{{\mathrm{GL}}}
    
    \newcommand{\Hom}{{\mathrm{Hom}}}
    
    \renewcommand{\Im}{{\mathrm{Im}}}
    \newcommand{\Ind}{{\mathrm{Ind}}}

    \newcommand{\ord}{{\mathrm{ord}}}

    \renewcommand{\mod}{\ \mathrm{mod}\ }

    \newcommand{\Res}{{\mathrm{Res}}}

    \newcommand{\SL}{{\mathrm{SL}}}

    \newcommand{\tr}{{\mathrm{tr}}}\newcommand{\tor}{{\mathrm{tor}}}
    
    \newcommand{\ur}{{\mathrm{ur}}}

    \newcommand{\ov}{\overline}

    \newcommand{\ra}{\rightarrow}

    \theoremstyle{plain}
    \newtheorem{thm}{Theorem}[section] 
    \newtheorem{lem}[thm]{Lemma}  \newtheorem{prop}[thm]{Proposition}

\theoremstyle{remark} \newtheorem{remark}[thm]{Remark}
\theoremstyle{remark} 
\theoremstyle{remark} 

    \newcommand{\cO}{\mathcal O}
    \numberwithin{equation}{section}
    
\begin{document}
\title{Dimension of the deformation space of ordinary representations in the cyclotomic limit
}

\author{Ashay Burungale, Laurent Clozel and Barry Mazur}
\address{Ashay A. Burungale: 
The University of Texas at Austin, Austin, TX 78712, USA.} 
\email{ashayburungale@gmail.com}
\address{Laurent Clozel: Math\'ematiques Universit\'e Paris-Saclay 91405 Orsay France} \email{laurent.clozel@math.u-psud.fr}

\address{Barry Mazur: Department of Mathematics, Harvard University, Cambridge, MA 02138,
USA}
\email{mazur@g.harvard.edu}
\maketitle
{\ \ \ \ \ \ \ \ \ \ \ \ \ \ \ \ \ \ \ \ \ \ \ \ \ \ \ \ \ \ \ \ \ \ \ \ \ \ \ \ \ \ \ \ \ \   \dedicatory{Dedicated to John H. Coates}}
\begin{abstract}
The weight two ordinary deformations are unobstructed in the cyclotomic limit under certain assumptions. 
We show that such an ordinary deformation ring over the cyclotomic tower can have arbitrarily large dimension. 
\end{abstract}

{
  \hypersetup{linkcolor=black}
  \tableofcontents
}

    \textit{ This work has evolved in the course of exchanges with Tony Feng, Dennis Gaitsgory and Michael Harris.  It is dedicated to the memory of John H. Coates (1945-2022), in appreciation of his achievements, his vision, and his multiple generosities.}
\section{Introduction}\label{Intr}
The paper concerns Iwasawa theory of ordinary deformation rings. 
\subsubsection*{Ordinary weight two (local) Galois representations}
    Let $k$ be a finite field,  $W(k)$ the Witt vectors of $k$ and $A$ a $W(k)$-algebra.
    
   For  $K$ a $p$-adic field with  $k $ its residue field, let $\ov{K}$ be an algebraic closure and $G_{K}= \Gal(\ov{K}/K )$.    A representation $\rho: G_{K} \rightarrow  \GL_{2}(A)$ is called \textit{ordinary of weight two} if it has the form

\begin{equation}\label{ordinary_def}  \left( \begin{array}{cc} \omega  \varepsilon & * \\ 0 & \varepsilon^{-1} \end{array} \right)\end{equation}
where $\varepsilon : G_{K} \rightarrow A^{\times}$ is unramified so that 
\begin{equation}\label{ds}\tag{ds}
\varepsilon^2 \neq 1, 
\end{equation}
and
$$\omega: G_{K} \rightarrow \mathbb{Z}_p^{\times} \rightarrow A^{\times}$$
is the cyclotomic character. 
\subsubsection*{(Global) Galois representations ordinary of weight $2$ at $p$} 
   Let $F$ be a totally real  number field. For an odd prime $p$, let $\zeta_p$ be a primitive $p$-th root of unity.   
   Let $F_{\infty}/F$ be the cyclotomic 
   $\BZ_p$-extension and $F_n$ its subextension so that $\Gal(F_{n}/F)\simeq \BZ/p^{n} \BZ$.  For a finite set of primes $S$ containing primes above $p$, let $F_S$ be the maximal extension of $F$ unramified outside $S$. 
   
    In a recent paper \cite{BC} Burungale and Clozel study the universal problem of deformations - in the cyclotomic limit $F_{\infty}/ F$ - of an ordinary, modular,  two dimensional representation $\bar{\rho}$ of $\Gal(F_S/F)$. 
    The aim of this paper is to consider the dimension of such deformation rings.

    We assume $F/\BQ$ unramified at $p$. Moreover, assume that $$\bar{\rho}: \Gal(F_S/F) \rightarrow \GL_{2}(k)$$ for $k$ finite, is an absolutely irreducible representation such that
\begin{itemize}
\item[(a)] $\bar{\rho}$ is ordinary of weight 2 at  the completions, $F_{\mathfrak p}$,  of $F$ at all primes ${\mathfrak p}$ dividing $p$.  
\item[(b)] The determinant is the cyclotomic character. 
\end{itemize}

Then there exists a universal deformation ring $R_n$ over $W(k)$, the \textit{ordinary deformation ring} for $F_n$, parametrizing ordinary liftings\footnote{as a representation of $\Gal(F_S/F_{n}).$} of weight 2 of $\bar{\rho}$ over algebras in the category $\widehat{\mathcal{C}}_{W(k)}$ of complete, local $W(k)$-algebras with residue field $k$.  See specifically Subsection 1.4  of \cite{BC}.  The ring $R_n$  is a complete Noetherian algebra (see Theorem 1.1 of \cite{BC}).
\subsubsection*{Deformations in the cyclotomic limit}
Define the deformation ring
 $$R_{\infty}=  \varprojlim R_n.$$ It belongs to $\widehat{\mathcal{C}}_{W(k)}$. It is not known to be Noetherian. (Compare \cite[pp. 354-357]{Hi}.) 

     The main result of \cite{BC} is that the deformation ring $R_{\infty}$ is formally smooth under some conditions. To recall, let $T \subset {\rm Ad}^{0}({\rho})$ denote an invariant lattice, where 
       $\rho$ is a lifting of $\bar{\rho}$ to characteristic $0$.

\begin{thm}\label{pre}
Assume $R_\infty$ is Noetherian.
   Moreover, assume that 
\begin{itemize}
\item[(i)] $\bar{\rho}$ is automorphic,  
\item[(ii)] $\bar{\rho}|_{G_{F(\zeta_{p})}}$ is adequate and 
\item[(iii)] $\mu(X^{1}(F,T^{*}(1))_{\tor})=0=\mu(X^{1}(F,T))_{\tor})$. 
\end{itemize}  

Then it is formally smooth, i.e.
$$R_\infty \cong W(k)[\![X_1,...,X_s]\!]$$
for some $s \geq 1$ (cf.~\cite[1.6]{BC}).
\end{thm}

We refer to \cite{Th} for the definition of `adequate' and to \cite{BC} for the modules $X^1$ and their $\mu$-invariants. Note that $R_\infty$ is known to be Noetherian if an adjoint $\mu$-invariant vanishes \cite[Cor.~5.11]{Hi}.
Recall that if $F=\BQ$, $\bar{\rho}$ is automorphic by the proof of Serre's conjecture \cite{KW}. 

In \cite{BC} it is also assumed that 
\begin{equation}\label{n-spl}
\text{The restriction of $\bar{\rho}$ to $ F_{\mathfrak{p}}$ is absolutely indecomposable for all $ \mathfrak{p} \vert p$.}
\end{equation}
As we explain in section \ref{s:uob}, this hypothesis is inessential. 

\subsubsection*{Dimension of the cyclotomic deformation space}
  The invariant $s$ appearing in Theorem \ref{pre} seems mysterious. 
  It satisfies $s\geq 1$. 
  
  A basic problem: to compute $s$ in given cases, and in particular  find examples with $s \geq 2$. In this note we will show that $s$ is arbitrarily large if arbitrarily large ramification is allowed. 
  An explicit example is presented in section \ref{s:ex} for the prime $p=5$. Then the case of general primes appears in section \ref{s:ex-g}. The latter relies on level raising, and pertinent deformation rings. 
  
   
   \vskip3mm

{\bf Acknowledgements}.
\noindent\vskip2mm
We thank George Boxer, Toby Gee, Haruzo Hida, Luc Illusie and Vincent Pilloni for helpful exchanges and suggestions. 
We also thank Christian Maire and Ariane Mézard  for instructive comments on a previous version. 
We are grateful to the referee for valuable suggestions. 
During the preparation of this paper, A.B. was partially supported by the NSF grants DMS-2303864 
and
DMS-2302064. 
\vskip5mm

\section{Ordinary deformations over cyclotomic tower}\label{s:uob} 
The aim of this section is to remove the hypothesis \eqref{n-spl} from the results of \cite{BC}. 

\subsection{ Setup} Let $p$ be an odd prime. Let $F$  be a totally real field of degree $d$ over $\mathbb{Q}$, unramified at $p$. All algebraic extensions of $F$ are contained in a fixed algebraic closure $\ov{F}$. 
Let $\zeta_p \in \ov{F}$ be a primitive $p$-th root of unity. 
  Let $F_{\infty}$ be the cyclotomic $\mathbb{Z}_p$-extension of $F$, and $F_n  \subset F_{\infty}$ the subextension of degree $p^n$ over $F$. 

Let $S$ be a finite set of places of $F$ containing the infinite and $p$-adic places.   
Let $F_S$ be the maximal extension of $F$ unramified outside $S$, and likewise define 
$F_{n,S}$. Put $\Gamma_0 = \Gal (F_S /F)$, $\Gamma_n = \Gal( F_{n,S} / F_n)$, and 
$\Pi=\Gal(F_{S}/F_{\infty})$. 
\subsubsection{Residual representation}\label{ss:res}
Let $k$ be a finite field of characteristic $p$. 
Let $\bar{\rho}: \Gamma_0 \rightarrow \GL_{2}(k)$
be an absolutely irreducible representation satisfying the following. 
\begin{itemize}
\item[(ord)] $\bar{\rho}$ is ordinary of weight 2.
\item[(irr$_{F(\zeta_{p})}$)] $\bar{\rho}|_{G_{F(\zeta_{p})}}$ is irreducible. 
\item[(det)] The determinant is the cyclotomic character.
\end{itemize}

Note that we do not assume the hypothesis (NS) of \cite{BC}, i.e. \eqref{n-spl} and so in particular allow residually CM cases satisfying the following: 
\begin{itemize}
\item[(rCM)] 
$$\bar{\rho}\simeq \Ind^{M}_{F}(\bar{\psi})$$
for $M/F$ a $p$-ordinary CM quadratic extension\footnote{i.e. any prime of $F$ above $p$ split in $M$} and $\bar{\psi}$ a finite order Hecke character over $M$. 
\item[(irr)] The character $\bar{\psi}/\bar{\psi}^{c}$ is non-trivial on $\Gal(\bar{\BQ}/F(\sqrt{(-1)^{\frac{p-1}{2}}p}))$ for $c\in\Gal(M/F)$ the non-trivial element and $\bar{\psi}^{c}:=\bar{\psi}\circ c$. 
\item[(ds)] The character $\bar{\psi}/\bar{\psi}^{c}$ is non-trivial on $\Gal(\bar{\BQ}_p/†M_{\mathfrak{p}})$ 
 for all $ \mathfrak{p}\ | \ p$.
\end{itemize}

\subsubsection{Ordinary deformations} 
Put $W=W(k)$. 
Write $\widehat{\mathcal{C}}_W$  for the category of complete local $W$-algebras with residue field $k$, and ${\mathcal{C}}_W$ for the subcategory of \textit{Artinian} objects in $\widehat{\mathcal{C}}_W$.

The ordinary deformation problem is well-posed in the cyclotomic tower. So 
for any non-negative integer $n$, there exists a universal deformation ring $R_n$ over $W(k)$, the \textit{ordinary deformation ring} for $F_n$ parametrising ordinary liftings (of weight 2) of $\bar{\rho}$ over algebras in  ${\mathcal{C}}_W$.   That is,  $R_n$ is the  complete local  $W(k)$-algebra in $\widehat{\mathcal{C}}_W$, defined up to unique isomorphism as representing the functor  

 $$A \ \mapsto\ {\rm  the \ set\ of\ ordinary\  liftings\  of\ } \bar{\rho} \ {\rm of \ weight\  2:}$$
\hskip160pt  $\xymatrix{\ &   \GL_2(A)\ar[d]\\
\Gamma_n\ar[r]^{\bar{\rho}}\ar[ur]^r  &\GL_{2}(k)}$
\vskip10pt
\noindent This ranges over $W(k)$-algebras $A$  in $\widehat{\mathcal{C}}_W$.
By ``lifting"  $r: \Gamma_n \to  \GL_2(A)$ we mean up to conjugation by $1+{\rm m}_AM_2(A)$ where ${\rm m}_A$ is the maximal ideal of $A$.

\vskip10pt
Note  that the natural homomorphism
$$R_{n+1}\to R_n$$
is surjective (cf.~\cite[Lem.~1.2]{BC}).
\vskip10pt 
 Put 
$$
R_{\infty}=\varprojlim R_{n}.
$$

We recall the following (cf.~\cite[Cor.~1.5]{BC}).

\begin{lem}\label{lim-univ}
$R_\infty$ represents the ordinary deformations of $\ov{\rho}|_{\Pi}$.
\end{lem}

\subsubsection{Main result}\label{subsub}

\begin{thm}\label{main}
Let $\bar{\rho}: G_F \rightarrow \GL_{2}(k)$
be an absolutely irreducible representation as above.
 Let $(\rho,V)$ be a deformation of $\bar{\rho}$ over the integer ring of a $p$-adic field, $V$ the underlying vector space and let $T \subset \Ad^0V$ be a $G_F$-stable lattice.
Assume $R_\infty$ is Noetherian.
   Assume further that 
   \begin{itemize}
  \item[(Aut)] $\rho$ is automorphic,   
  \item[(ad$_{F(\zeta_{p})}$)] $\bar{\rho}|_{G_{F(\zeta_{p})}}$ is adequate and 
  \item[($\mu$)] $\mu(X^{1}(F,T^{*}(1))_{\tor})=0=\mu(X^{1}(F,T)_{\tor})$ 
   \end{itemize} 

Then it is formally smooth, i.e.
$$R_\infty \simeq W(k)[\![X_1,...,X_s]\!]$$
for some $s \geq 1$.
\end{thm}

We refer to \cite[\S4]{BC} for the definition of Selmer groups in the hypothesis ($\mu$). 

The result is proved by the same method as in \cite{BC}. We only indicate the differences. Since 
$\ov{\rho}$ is allowed to split at primes above $p$, the essential difference is calculation of local cohomology at such primes.

   \subsection{Local cohomology}\label{s:local}

In this section\ $K=F_{n,\mathfrak{p}}$ is local and $\bar{\rho}$ is an ordinary representation of $G_K$ as above. 
For simplicity we write $d$ for $d_{\mathfrak{p}} = [F_{\mathfrak{p}}:\BQ_p].$

Since the locally indecomposable case is covered by \cite[\S2]{BC}, we assume that 
the restriction of $\ov{\rho}$ to the decomposition group at $\fp$ is of the form 
\begin{equation}\label{ordinary_def}  \left( \begin{array}{cc} \bar{\omega} \bar{\varepsilon} & 0 \\ 0 & \bar{\varepsilon}^{-1} \end{array} \right)\end{equation}
where $\ov{\varepsilon} : G_K \rightarrow k^{\times}$ is unramified with $\bar{\varepsilon}^2 \neq 1$ 
and
$\bar{\omega}: G_K \rightarrow \BZ_p^{\times} \rightarrow \BF_{p}^{\times}$
the mod $p$ cyclotomic character.

Let $V=k^2$ be the space of $\bar{\rho}$, and $W = \End^0(V)$ be the space of traceless endomorphisms of $V$. It is endowed with the natural representation $\Ad^0\bar{\rho}$. Let $(\Ad^0\bar{\rho})(1)$  be the Tate twist.

The main result is Proposition \ref{spld} below. 

\subsubsection{Local cohomology of the adjoint}
\vskip2mm
\subsubsection*{Filtration}
Let  $W_0 \subset W_1 \subset W_2=W$ be the filtration of $W$:
\begin{equation}\label{udl}
W_{0}=\bigg{\{}\left( \begin{array}{cc} 0 & * \\ 0 & 0 \end{array}\right) \bigg{\}}, 
\qquad 
W_{1}=\bigg{\{}\left( \begin{array}{cc} * & * \\ 0 & * \end{array} \right) \bigg{\}}
\end{equation}
preserved by $G_K$. Then as $G_K$-modules,
$$W_0  \cong k[\varepsilon^2 \omega], \qquad W_1/W_0 \cong k[1], \qquad W_2/W_1 \cong k[\omega^{-1} \varepsilon^{-2}] $$ 
for $1$ being the trivial character. 

The exact sequence
\begin{equation}\label{fil_1}
0 \rightarrow W_0 \rightarrow W_1 \rightarrow W_1/W_0 \rightarrow 0
\end{equation}
induces
\vskip2mm
$H^0(K, W_1) \rightarrow H^0(K, W_1/W_0) \rightarrow H^1(K, W_0) \rightarrow H^1(K, W_1) \rightarrow H^1(K, W_1/W_0) \rightarrow H^2(K, W_0)  \rightarrow H^2(K, W_1)  \rightarrow H^2(K, W_1/W_0) \rightarrow 0.$
\vskip2mm
Write $h^i(K, -) = \dim_{k}H^i(K,-).$
\begin{lem}\label{dim}
\item[(i)] $h^0(K,W_1/W_0)= 1$, 
\item[(ii)] $h^1(K,W_0)= p^n d$, 
\item[(iii)] $h^1(K, W_1) = 2 p^n d +1$, 
\item[(iv)] $h^1(K, W_1/W_0) = p^n d +1$, 
 \item[(v)] $ h^2(K,W_1)=  0 $.

 \end{lem}
 \begin{proof}
 Write $V'= V^*(1)$. Then $$W_0' \cong k[\varepsilon^{-2}], \qquad (W_1/W_0)' \cong k[\omega].$$ By Tate duality we see that $h^2(K, W_0) = h^2(K, W_1/W_0)= 0 .$ This implies (v). 
 
 Note that (i) is obvious since the extension is split. 
 The map $H^1(K, W_1) \rightarrow H^1(K, W_1/W_0)$ is surjective since $H^2(K,W_0)=0$.  Now the formulas (ii)-(iii) follow from Tate's Euler-Poincar\'{e} formula and (iv) from the exact sequence. 
 \end{proof}
 
For $U$ a representation of $G_K$ on a $k$-vector space, put 
 $$H^1_{\mathrm{nr}}(K,U) = \ker \{H^1(G_K,U) \rightarrow H^1(I_K,U)\}$$
where $I_K$ is the inertia. We \textit{define} the unramified classes $H^1_{\ur}(K, W_1)$ to be the inverse image of $H^1_{\mathrm{nr}}(K, W_1/W_0).$ 
  
  We have the exact sequence\footnote{Note that \eqref{fil_1} is a split exact sequence by the hypothesis \eqref{ordinary_def}.}
 \begin{equation}\label{udses}
 0 \rightarrow H^1(K, W_0)   \rightarrow H^1(K, W_1)    \rightarrow H^1(K, W_1/ W_0) \rightarrow  0 
\end{equation}
where the corresponding $k$-dimensions are $ ( p^n d, 2 p^nd+1, p^n d+1)$. Since $W_1/W_0$ is with trivial $G_K$-action, $H^1_{\rm{nr}}(K, W_1/W_0) \cong k$. Thus
$$
\dim_{k}H^1_{\ur}(K, W_1)= p^n d + 1.
$$

 Now the exact sequence
 \begin{equation}\label{fil_2}
 0 \rightarrow W_1 \rightarrow W \rightarrow W/W_1 \rightarrow  0 
 \end{equation}
 induces
 \begin{equation}\label{dlses}
  0 \rightarrow H^1(K, W_1)  \rightarrow H^1(K,W) \rightarrow H^1(K, W/W_1) \rightarrow 0
 \end{equation}
 by Lemma \ref{dim}. 
   
   We \textit{define} $H^1_{\ord}(K, \Ad^0 \bar{\rho})$ as the image of $H^1_{\ur}(K, W_1)$ in $H^1(K,W)$. We also note the vanishing of $H^2(K, \Ad^0 \bar{\rho})$ by the analogue of (\ref{dlses}) for $H^2$, and Tate duality for $W/W_1$.
   
   We summarise the results obtained so far:
   \begin{lem}\label{lc1}
   \item[(i)] $\dim_{k} H^0(K, \Ad^0 \bar{\rho}) = 1$ and $H^2(K, \Ad^0 \bar{\rho}) = 0.$
   \item[(ii)] $\dim_{k} H^1_{\ord}(K, \Ad^0 \bar{\rho}) = p^n d + 1$.
   \item[(iii)] $\dim_{k} H^1(K, \Ad^0 \bar{\rho}) = 3 p^n d + 1.$
 
   \end{lem}

   Now consider the extension $K=F_{n, \mathfrak{p}} = K_n$ of $K_0= F_{\mathfrak{p}}$, whence an action of $\Delta_n= \Gal(K_n/K_0)$ on the cohomology groups $H^*(K_n, -).$
\begin{lem}\label{sec:freeness} 
As $k[\Delta_{n}]$-modules, we have
$$H^1(K_n, (\Ad^0\bar{\rho}) (1)) \simeq k[\Delta_{n}]^{3d}\oplus k$$
with $k$ being the trivial $k[\Delta_{n}]$-module. 
\end{lem}
\begin{proof}   
   By Tate's Euler characteristic formula, $\dim_{k} H^1(K_n, (\Ad^0\bar{\rho}) (1)) = 3 p^n d + 1$.

The lemma is a consequence of the exact sequences \eqref{fil_1} and \eqref{fil_2} being split. 
Indeed, as $k[\Delta_{n}]$-modules, we have
$$
\text{
 $H^{1}(K_{n},W_{1}/W_{0}(1))\simeq k[\Delta_{n}]^{d}\oplus k$, and  
$H^{1}(K_{n},W_{0}(1))\simeq H^{1}(K_{n},W_{2}/W_{1}(1))\simeq k[\Delta_{n}]^d$.
}
$$
The latter follows by considering the space of $\Delta_n$-coinvariants, and the corresponding dimensions. Likewise, the former holds by the split exact sequence \eqref{fil_1}.

 \end{proof}  
   
   We now consider the subspace $H^1_{\ord}(K_n, \Ad^0 \bar{\rho})$, of dimension $p^n d + 1$. Note that the filtration $W_i$ of $W$ gives rise to cohomology  spaces on which $\Delta_n$ acts.
   
   \begin{lem}\label{inv}
   
  $H^1_{\mathrm{nr}}(K_n, W_1/W_0)$, $H^1_{\ur}(K_n, W_1)$ and $H^1_{\ord}(K_n, \Ad^0\bar{\rho})$ are invariant under the action of $\Delta_n$.
   
   \end{lem}

 \begin{lem}\label{spl}
 $H^1_{\ord}(K_n, \Ad^0 \bar{\rho})$ is isomorphic, as a $k[\Delta_n]$-module, to  $$k[\Delta_n]^d  \oplus k$$ with $k$ being the trivial $k[\Delta_n]$-module.
  \end{lem}
 \begin{proof}  
    The exact sequence (\ref{udses}) yields 
    \begin{equation}\label{urses}
    0 \rightarrow H^1(K_n, W_0) \rightarrow H^1_{\ur}(K_n, W_1) \rightarrow H^1_{\mathrm{nr}}(K_n, W_1/W_0) \rightarrow 0
   \end{equation}
 with $H^1_{\ur}(K_n, W_1) \cong H^1_{\ord}(K_n, \Ad^0\bar{\rho})$ and the dimensions being $(p^nd, p^n d + 1, 1)$. The argument given for Lemma \ref{sec:freeness} shows that $H^1(K_n, W_0)$ is free of rank $d$ over $k[\Delta_n]$. Then the assertion follows from the proof of \cite[Lem.~2.5]{BC}.     
\end{proof}

 Consider now $\pi: \Delta_{n+1} \twoheadrightarrow \Delta_n$. This induces a natural map $k[\Delta_n] \hookrightarrow k[\Delta_{n+1}]$, $f(\delta) \mapsto f(\pi \delta),$ dual to the projection of Iwasawa theory. It is equivariant under the action of $\Delta_{n+1}$, acting on $k[\Delta_n]$ via the quotient map.
 
 \begin{lem}\label{ocp}
 The restriction $H^1_{\ord}(K_n, \Ad^0 \bar{\rho}) \rightarrow H^1_{\ord}(K_{n+1}, \Ad^0 \bar{\rho})$ is injective. It is compatible with the splitting of Lemma \ref{spl}, and equivariant for the action of $\Delta_{n+1}$. 
 \end{lem}
This is the content of \cite[Lem.~2.6]{BC}.
 
 \subsubsection*{Local cohomology, dualised}
 
 \vskip2mm
 
 We now use the Tate pairing
  $$ H^1(K_n, \Ad^0\bar{\rho}) \times H^1(K_n, \Ad^0\bar{\rho}(1)) \rightarrow k.$$

   Let $ H^1_{\ord, \bot} \subset H^1(K_n, \Ad^0\bar{\rho}(1))$ be the orthogonal space of $H^1_{\ord}$. We set
$$ H^1_{\ord, *}(K_n, \Ad^0\bar{\rho}(1)) =  H^1(K_n, \Ad^0\bar{\rho}(1))/ H^1_{\ord,\bot}.$$
      So this is naturally dual to $H^1_{\ord}$. When $K_n$ is concerned, we write $H^1_{\ord , n}$ etc. We can take the limit of these spaces under corestriction. In fact we obtain naturally a diagram
      $$
\xymatrix{
 H^1_{\ord, *, n+1} \ar[d]\ar[r]^{\cong }&(H^1_{\ord, n+1})^* \ar[d]\\
 H^1_{\ord, *, n}\ar[r]^{\cong } &(H^1_{\ord, n})^*\\
}
$$
 where the surjection on the right comes from the previous injection (Lemma \ref{ocp}) and the surjection on the left completes the diagram. 
Note that for $\beta \in H^1_{\ord, *, n+1}$ and $\alpha \in H^1_{\ord,n}$, 
 $$
 ( \Cor ~ \beta, \alpha ) = ( \beta, \Res~ \alpha ) 
 $$ 
(cf.~\cite[\S2.2]{BC}).

 We now dualise the expression of $H^1_{\ord}(K_n, \Ad^0\bar{\rho})$ obtained in Lemma \ref{spl}. 
We have
$$H^1_{\ord}(K_n, \Ad^0\bar{\rho})^* \cong (k[\Delta_{n}]^{d} \oplus k)^* \cong  (k[\Delta_{n}]^{d})^* \oplus k.$$
 If we restrict to $K_{n+1}$, the corresponding map $k \rightarrow k$ is an isomorphism as was seen in the proof of Lemma \ref{spl}. In view of the preceding paragraph we deduce the following.

  \begin{prop}\label{spld} Put
  $$
\Omega=\varprojlim k[\Gamma_{n}] \simeq k[\![T]\!].
$$  
  Then, as an $\Omega$-module, we have
     $$\varprojlim H^1_{\ord,*,n}(K_n, \Ad^0\bar{\rho}(1)) \cong \Omega^d \oplus k .$$
  \end{prop}
  
    \subsection{Main result}\label{s:ub} 
    In this section we consider the vanishing of the second ordinary global Galois cohomology for adjoint over the cyclotomic tower. 
    
        Let the notation and hypotheses be as before.  
        \begin{prop} \label{free} 
         Suppose that 
   \begin{itemize}
  \item[(Aut)] $\bar{\rho}$ is automorphic,   
  \item[(ad$_{F(\zeta_{p})}$)] $\bar{\rho}|_{G_{F(\zeta_{p})}}$ is adequate (\cite{Th}) and 
  \item[($\mu$)] $\mu(X^{1}(F,T^{*}(1))_{\tor})=0=\mu(X^{1}(F,T)_{\tor})$ for $T$ a lattice in $\Ad^{0}(\rho)$ with $\rho$ arising from an automorphic lift. 
   \end{itemize} 
         
        Then, $H^{1}_{\mathrm{ct}}(\Gamma_{0},\Ad^{0}\bar{\rho}(1) \otimes \Omega)$ is free over $\Omega$ of rank $[F:\BQ]$.
        \end{prop}
        \begin{proof}  We first show that $H^{1}_{\mathrm{ct}}(\Gamma_{0},\Ad^{0}\bar{\rho}(1) \otimes \Omega)$ is free as an $\Omega$-module: it suffices to show that the corresponding $\Omega$-torsion submodule 
        \begin{equation}\label{tor}
         (\Ad^{0}\ov{\rho}(1)\otimes_{k} \Omega)^{G_F} 
         \end{equation}
         vanishes (cf.~\cite[p.~12,~Lemma]{PR}). This is a consequence of the hypothesis that $\ov{\rho}$ is an irreducible $k[G_F]$-module. Accordingly, either ${\rm Ad}^{0}\ov{\rho}$ is an irreducible $k[G_F]$-module or $\ov{\rho}$ is dihedral. 
         \begin{itemize}
         \item First, suppose that ${\rm Ad}^{0}\ov{\rho}$ is an irreducible $k[G_F]$-module. Note that \eqref{tor} is identified\footnote{This identification is independent of irreducibility.} with the set of all $k[G_{F}]$-homomorphisms $\Hom_{k}(\Omega, k) \ra {\rm Ad}^{0}\ov{\rho}(1)$, where $G_F$ acts on $\Omega$ via the surjection $G_{F}\twoheadrightarrow \Gal(F_{\infty}/F)$.
Since the action of $G_F$ on $\Hom_{k}(\Omega, k)$ is abelian and the representation
of $G_F$ on ${\rm Ad}^{0}\ov{\rho}$ is irreducible and non-abelian, such a non-trivial
homomorphism does not exist.         
         \item In the dihedral case, \eqref{tor} clearly vanishes. 
         \end{itemize}
        
        To calculate the rank, the proof of \cite[Thm.~5.2]{BC} applies.

        \end{proof}

        The main theorem:
 \begin{thm} \label{unob}
   Suppose that 
   \begin{itemize}
  \item[(Aut)] $\bar{\rho}$ is automorphic,   
  \item[(ad$_{F(\zeta_{p})}$)] $\bar{\rho}|_{G_{F(\zeta_{p})}}$ is adequate (\cite{Th}) and 
  \item[($\mu$)] $\mu(X^{1}(F,T^{*}(1))_{\tor})=0=\mu(X^{1}(F,T)_{\tor})$ for $T$ arising from an automorphic lift. 
   \end{itemize} 
 Moreover, suppose that  $R_{\infty}$ is Noetherian.   Then  $\varinjlim H^2_{\ord}(\Gamma_n, \Ad^{0}\bar{\rho}) = 0$;
             in particular $$R_{\infty} \simeq W(k)[\![X_1, \dots, X_s]\!].$$ 
   \end{thm}
    \begin{proof}  
    The proof of \cite[Thm.~5.2]{BC} is based on Lemma 2.7 and Proposition 5.1 of {\it{loc. cit.}}, which are proved under the hypotheses (NS).  As seen in Propositions \ref{spld} and  \ref{free},     
    the assertions also hold in our setting, and so the proof of \cite[Thm.~5.2]{BC} applies.

    \end{proof}

\section{An explicit example for $p=5$}\label{s:ex}

\textit{In the rest of the paper, we assume that $\bar{\rho}$ satisfies the conditions in Theorem~\ref{pre}.}

\subsection{The example} 
Our base field is $\BQ$. We choose $p=5$. 

\subsubsection{Setting}\label{ss:st}
 The curve $\CX= \CX$(5) - the modular curve of full level $\Gamma(5)$ - has genus $0$. The elliptic  curve $E$ with label $37a1$ in \cite{Cr}   (and labelled $37.a1$ in \cite{LMFDB}) is $5$-ordinary, of level $37$. It defines $\bar{\rho}$. 
 
 Note that the $p$-distinguished hypothesis \eqref{ds} holds since $\varepsilon(\Fr_{p})^{-1}$ is the $p$-unit root of the Hecke polynomial $X^2+2X+5$. According to \cite{LMFDB}, $\Im(\bar{\rho}) = \GL_{2}(\BF_5)$. 
Thus $\Im\, (\bar{\rho})|_{\Gal(\ov{\BQ}/\BQ(\zeta_5))}=\SL_2(\BF_5)$. In particular  this restricted representation is {\it big} in the sense of \cite{CHT} (Cor.~2.5.4 of this paper) and therefore adequate \cite{Th}. We assume that the condition (iii) in Theorem~\ref{pre} is satisfied.

    Let $\CX(\bar{\rho}) = \CX'$ be the twist of $\CX$ by $\bar{\rho}$ (cf.~\cite[p.~543]{Wi}). It has one rational point, and therefore 
    infinitely 
    many rational points yielding elliptic curves over $\BQ$ with the same Galois representation mod  5; i.e., $$\bar{\rho}: \Gal({\ov{\bf Q}}/ {\bf Q}) \to \GL_{2}(\BF_5).$$ Consequently:
  \vskip10pt \begin{lem}\label{lem1} There are infinitely many elliptic curves over $\BQ$  having $\bar{\rho}$ as mod $5$ Galois representation that have different conductors (and hence are  mutually non-isogenous). \end{lem}
  \begin{proof} Given the discussion above, and since the natural mapping to the $j$-line is of finite degree,  the image ${\mathcal E}\subset {\bf P}^1({\bf Q})$ of the infinite set $\CX'(\BQ)$ to the $j$-line is infinite. Since for any conductor $N$ there are only finitely many mutually non isomorphic elliptic curves over ${\bf Q}$ of the same conductor, any elliptic curve over $\BQ$ has only a finite number of isogenous $\BQ$-elliptic curves, imposing the equivalence relation $E\sim E'$  on the elements $E,E' \subset {\mathcal E}$ if $E$ and $E'$ have the same conductor. Therefore, we have infinitely many different equivalence classes in  ${\mathcal E}$.\end{proof} 
  \vskip10pt
  \begin{prop}\label{prop:seq}  There is an infinite sequence \begin{equation}\label{infell}E=E_1, E_2,\dots, E_m,\dots\end{equation} of  elliptic curves over $\BQ$, and  {\it distinct} prime numbers $p_1,p_2,\dots, p_i,\dots$ different from $p=5$ with the following properties.
    If $\rho_i:\Gal({\ov {\bf Q}}/ {\bf Q}) \to {\rm Aut}\big\{T_5(E_i)\} \simeq \GL_2({\BZ}_5)$ is the Galois representation given by the action of $\Gal({\bar {\bf Q}}/ {\bf Q})$ on the  $5$-adic Tate module  built on $5$-power torsion points of $E_i$ then:
\begin{enumerate}\item The residual representation  ${\bar \rho}_i: \Gal({\ov {\bf Q}}/ {\bf Q}) \to {\rm Aut}_{\BF_{5}}(E_i[5])$  is isomorphic to ${\bar\rho}$. 
  \item  The restriction of  $\rho_i$ to $\Gal({\ov {\bf Q}}_{p_i}/ {\bf Q}_{p_i}) $  is ramified, while the restriction of  $\rho_j$ to $\Gal({\ov {\bf Q}}_{p_i}/ {\bf Q}_{p_i}) \to  \GL_2({\BZ}_5)$ is unramified for $j <i$.
  \end{enumerate}\end{prop}
 \vskip10pt
 \begin{proof}  Begin by forming  an infinite sequence  of elliptic curves  as given by Lemma \ref{lem1}  ordered so that they have strictly increasing conductors: $N_1 < N_2 <\dots < N_n<\dots$.   This sequence satisfies condition (1) of our proposition. 
 
 Now winnow it as follows. Supposing that we have already chosen $E_1,E_2,\dots, E_m$  so that for every $j < i \le m$   the restriction of $\rho_i$ to $\Gal({\ov {\bf Q}}_{p_i}/ {\bf Q}_{p_i})$ is ramified while the restriction of  $\rho_j$ to $\Gal({\ov {\bf Q}}_{p_i}/ {\bf Q}_{p_i})$ is unramified  (i.e., suppose that condition (2) holds up to $i=m$).  
 
 \vskip10pt
     
      {\bf To prove:}  there is an index $n$ such and a prime $p_n$ dividing $N_n$ such that $p_n$ does not divide $\prod_{i=1}^mN_i$. \\

       The reason for this is that if there is no such prime, all the conductors $\{N_i; i = 1,2,\dots\}$ would be  expressible as products of powers of a finite set of primes--but this is not possible given that there is a finite maximum   exponent for any prime power dividing any conductor.

        Now, casting out all the elliptic curves $E_j$ with indices $i$ that lie properly between $m$ and $n$, and relabelling appropriately---i.e., the formerly named $E_n, N_n$ and $p_n$ are relabeled as  $E_{m+1}, N_{m+1}$ and $p_{m+1}$---we have (by the N{\'e}ron-Ogg-Shafarevich criterion) that 
 for $i \le m$ the representation $\rho_i$ is unramified when restricted to   $\Gal({\ov {\bf Q}}_{p_{m+1}}/ {\bf Q}_{p_{m+1}})$ since $p_{m+1}$ does not divide the conductor $N_i$ of $E_i$   (for $i \le m$).  By the same N{\'e}ron-Ogg-Shafarevich criterion,   $\rho_{m+1}$ is ramified when restricted to   $\Gal({\ov {\bf Q}}_{p_{m+1}}/ {\bf Q}_{p_{m+1}})$.
 \end{proof}
 We note that for all these curves, the image of $\Gal(\bar{\BQ}/\BQ)$ in $\Aut(T_{5} E)$ equals $\GL_2(\BZ_5)$ by a theorem of Serre (cf.~\cite[IV.~3.4]{Se}) 
 since the residual image is full. \vskip10pt 
 {\it Side comment:}
     In contrast with the case we are considering (i.e., $p=5$) when $p=3$ the $p$-distinguished hypothesis \eqref{ds} cannot hold. Nevertheless, note  that  $\frak{X} (3)$ has genus 0 and we have examples, such as the curve $E=11.a.1$ in \cite{Cr},  that satisfy all the same properties  (relative to $p=3$) that   $37a1$ satisfies (relative to $p=5$) save for \eqref{ds}.

 \vskip10pt

Now let $S_0=\{p,p_{0}\}=\{5,37\}$. 
 Using Proposition~\ref{prop:seq} and extracting a subsequence if necessary, we can find an infinite sequence $(p_m)$ and an increasing, infinite sequence $(S_m)$ of finite sets of primes containing $S_0$ and $p_m$, with the following property.

  \begin {thm}  There is an infinite sequence of pairs defined over ${\BQ}$,  $$\{\tilde{E}_m{_{/{\bf Q}}};\  \alpha:\  E[p]_{/{\bf Q}}\stackrel{\simeq}{\to} \tilde{E}_m[p]_{/{\bf Q}}\}_{m\in\BZ_{\geq 0}}$$   
   with $E_0=E$ as above; and an infinite set of (distinct) primes $p_m \ne p$ such that starting with $S_0$ as above, and defining inductively $S_m:= S_{m-1}\sqcup\{p_m\}$ we have that $\tilde{E}_m$ has good reduction at primes  outside $S_m$ and bad reduction at  $p_m$.   \end{thm}

Recall that $\bar{\rho}$ is ordinary, i.e, at the prime $p=5$: 

\begin{equation}\label{eq:res}
\bar{\rho} \simeq \begin{pmatrix}
\omega \varepsilon &*\\
&\varepsilon^{-1}
\end{pmatrix},\ \varepsilon\ \mathrm{unramified,}\ \varepsilon^2\not= 1
\end{equation}
and that this notion is defined if $\rho$ is a representation on a $W(\BF_5)$--algebra.

  For each of the curves $\tilde{E}_m$, the representation of $\tilde{E}_m[5]$ is ordinary by construction. However, 
    this does not\footnote{In fact, this is true if $p>7$; however George Boxer found counterexamples for $p\leq 7$. The argument is quite interesting and we hope that a proof will be published.} imply that so is the representation on $T_5(\tilde{E}_{m})$.

\subsubsection{Variant}    We now refine the foregoing construction so as to guarantee $5$-ordinarity. This wonderfully simple argument was indicated to us by Vincent Pilloni.
    
    To begin, note that $5$-ordinariness is an open condition in the $5$-adic topology in $X(\ov{\rho})$. Indeed, this holds since the projection $X(\ov{\rho}) \rightarrow X(1)$  is a finite morphism and  $5$-ordinariness is  visibly an open condition in the $5$-adic topology of $X(1)$.    
    
    So a neighborhood $U \subset X(\ov{\rho})(\BQ_5)$ of $E_0$ contains only ordinary points defined over $\BQ_5$.  On the other hand, $U \cap X(\ov{\rho})(\BQ)$ is infinite since $X(\ov{\rho})$ is isomorphic to $\BP^1$ over $\BQ$.  
        Pick its infinite subset $\{E_m\}_{m \in \BN}$ with $\{\rho_m\}$ the associated $5$-adic Galois representations. We can now use the construction in subsection \ref{ss:st} for the points of this set.  With this condition we now have:    
    
  \begin{prop}  
    For each curve $E_m$, the representation $\rho_{m}$ of $G_\BQ$ on $T_5 (E_m) \cong \BZ_5^2$ is ordinary and unramified outside~$S_m$.
    \end{prop}

\begin{prop}\label{prop:ram}
Let $G_m=\Gal(\BQ_{S_m}/ \BQ_\infty)$. 
Then $\rho_m|_{G_m}$ is ramified at $\fp_m^{(\infty)}$, i.e. $\rho_m|_{I_m^\infty}$ is non--trivial.
\end{prop}
\begin{proof}
Recall that $\rho_m\vert_ {G_m}$ is ordinary. Consider a sequence of primes $$(p_m) \subset \fp_m^{(1)} \subset \fp_m^{(2)} \subset\cdots \fp_m^{(\infty)}$$ in the tower $\BQ\subset \BQ_1 \subset \cdots \subset \BQ_\infty$. Since $p_m$ is unramified in $\BQ_\infty$, we have isomorphisms $I_m^0 = I_m^1=\cdots = I_m^\infty$ where $I_m^r \subset \Gal(\BQ_{S_m}/\BQ_r)$ is the inertia group for $\fp_m^{(r)}$. The assertion follows. 
\end{proof}

\subsection{Growth of deformation rings} We analyse deformation rings arising in the context of  preceding subsection.

To begin, in this situation the conditions of  Theorem~\ref{pre} are satisfied.
We then have, for all $m$, the deformation ring associated to $\bar{\rho}$ and $S_m$, denoted by $R_\infty^m$. It parametrises ordinary representations of $\Gal(\ov{\BQ}/\BQ_\infty)$ lifting $\bar{\rho}
$ and unramified outside $S_m$. In particular, for $m'\le m$, we get a natural map
\begin{equation}\label{eq:can}
{\rm can} : R_\infty^m \rightarrow R_\infty^{m'}.
\end{equation}

By Theorem~\ref{pre}, $R_\infty^m$ is, for all $m$, isomorphic to $\BZ_5[\![X_1,\ldots , X_{s_m}]\!]$ with $s_m\ge 1$, if it is Noetherian.

\begin{prop}
The homomorphism \eqref{eq:can} is surjective.
\end{prop}
\begin{proof} Write $\CR$ for $R_\infty^m$, $\CR'$ for $R_\infty^{m'}$ and 
$G$ for $\Gal(\ov{\BQ}/\BQ_\infty)$. 

We then have maps
\begin{equation}\label{eq:com}
G \stackrel{{\rm univ}}{\to} \GL_{2}(\CR) \stackrel{{\rm can}}{\to} \GL_{2}(\CR')
\end{equation}
with composition ${\rm univ'}$, where 
${\rm univ, univ' }$ are the universal representations. Let $\CR_0\subset \CR'$ be the image of \eqref{eq:can} and $j:\CR_0\rightarrow \CR'$ the injection. We rewrite \eqref{eq:com}~as
$$
\xymatrix{
G \ar@{=}[d]\ar[r]^{{\rm univ}\qquad}&\GL_2(\CR)\ar[d]^{{\rm can}}\\
G \ar[r]^{{\rm univ'}\qquad}&\GL_2(\CR');
}
$$
since $\CR\rightarrow \CR'$ quotients through $\CR_0$ this yields
$$
\xymatrix{
    G_m \ar[d]\ar[r]\ar[rd]&\GL_2(\CR_0)\ar[d]\\
    G_{m'} \ar[r]&\GL_2(\CR'),
}
$$
taking account of the ramification of ${\rm univ,univ'}$. The diagonal map factors through $G_{m'}$, so it is unramified outside $S_{m'}$. The inclusion $\GL_2(\CR_0)\rightarrow \GL_2(\CR')$ is injective; this implies that the natural map $G_m\rightarrow \GL_2(\CR_0)$ is unramified outside $S_{m'}$. Therefore it induces a canonical map $\sigma:\CR'\rightarrow \CR_0$ such~that
$$
\xymatrix{
    G\ar[r]\ar[rd]^{{\rm univ'}}&\GL_2(\CR_0)\\
    &\GL_2(\CR')\ar[u]_{\sigma}
}
$$
commutes; finally we have a diagram
$$
\xymatrix{
G \ar@{=}[d]\ar[r]&\GL_2(\CR_0)\ar[d]^{j}\ar@{=}[r]&\GL_2(\CR_0)\\
G \ar[r]&\GL_2(\CR')\ar@{=}[r]&\GL_2(\CR')\ar[u]_{\sigma}.
}
$$
In particular the diagram 
$$
\xymatrix{
G \ar@{=}[d]\ar[r]^{{\rm univ'}\qquad}&\GL_2(\CR')\ar[d]^{j\circ \sigma}\\
G \ar[r]^{{\rm univ'}\qquad}&\GL_2(\CR')
}
$$
also commutes, which implies that $j\circ \sigma={\rm id}_{\CR'}$, hence that $j$, and in turn $\CR \rightarrow \CR'$, is surjective.
\end{proof}
\begin{lem}\label{lm:nsj}
The homomorphism \eqref{eq:can} is not bijective.
\end{lem}
\begin{proof}
Assume it is bijective. Then the universal representation ${\rm univ}$ factors through $\CR'$, and is  therefore unramified outside $S_{m'}$. In particular $\rho_m\mid _{G_m}$ would be so, but this contradicts Prop.~\ref{prop:ram}.
\end{proof}

The main result of this subsection is the following. 
\begin{thm}\label{thm:gr}
Assume deformation rings $R_\infty^m$ are Noetherian for any $m\geq 0$.
Then $R_\infty^m$ is formally smooth of dimension at least $m+2$.
\end{thm}
\begin{proof}

The natural map $R_\infty^m {\overset{\varphi}{\rightarrow}} R_\infty^{m'}$ is surjective, not an isomorphism by Lemma~\ref{lm:nsj}. It follows that $$s_m\ge s_{m'}:= s.$$ Suppose that $s_m=s_{m'}$, both rings being isomorphic to $\BZ_5[\![X_1,\ldots X_s]\!]$ by Theorem~\ref{pre}. 

Let $I_r$ be the kernel of the $r$--th iterate of $\varphi$. Since $I_1\not=0$ and $\varphi$ is onto, $I_2$ contains an element mapping (by $\varphi$) to a non--zero element of $I_1$. In  particular $I_2$ properly contains $I_1$. The same argument applies to $\varphi^r : I_{r+1}\rightarrow I_1$, showing that $I_{r+1}$ properly contains $I_r$. This is impossible since the rings are Noetherian.

Here is an alternate argument due to the referee: Since $\ker \varphi \neq 0$ and $R_{\infty}^m$ is a domain, $R_{\infty}^m/\ker\varphi$ has Krull dimension at most $s_m$. On the other hand, $R_\infty^{m'}$ has Krull dimension $1+s_{m'}$.
 The contradiction concludes the proof. 
\end{proof}
\begin{remark}
The Noetherian-ness of $R_\infty^{0}$ may be checked by verifying the $\mu$-invariant hypothesis in \cite[Cor.~5.10]{Hi}. Then $R=T$ results seem to imply that $R_\infty^m$ is Noetherian for any $m$. 
\end{remark}

\section{General $\bar{\rho}$}\label{s:ex-g}

We now consider the general case of an odd prime $p$ and residual Galois representation $\bar{\rho}$ satisfying the conditions of Theorem~\ref{pre} for $F=\BQ$.

Assume that $\bar{\rho}$ is odd and irreducible. By the proof of Serre's conjecture \cite{KW}, 
$\bar{\rho}$ is modular. 
So it is associated to a newform $f$ on $\Gamma_1(N)$ of weight $\kappa$, where $(N,\kappa)$ are given by Serre's recipe. In particular $N$ is coprime to $p$. By \cite[\S4.1]{Edix}, we have $\kappa\in\{2,p+1\}$.

\subsection{The case $\kappa=2$}\label{ss:wt-2}
 We assume that the representation of $G_\BQ$ on the space $T\cong \CO_L^2$, where $L/\BQ_p$ is a field of coefficients for $f$, verifies (iii) of Theorem~\ref{pre}; we also assume that $\bar{\rho}|_{\BQ(\zeta_{p})}$ is adequate. For simplicity (cf.~\cite[Thm.~3]{DT}), assume that $p\geq 5$.  
 
 Consider the primes $\ell$, not dividing $Np$, verifying the condition: 
 \begin{equation}\label{R}\tag{R}
 \text{The characteristic polynomial  of $\bar{\rho}(\Frob_\ell)$  is of the form $X^2 + aX+b$, with $\ell a^{2}=(\ell+1)^{2}b$.}
 \end{equation}
Here $\Frob_\ell$ is the arithmetic  Frobenius. 
The condition is equivalent to $\bar{\rho}(\Frob_\ell)$ having two eigenvalues $\alpha,\beta$ with $\alpha=\ell \beta$. We call these \textit{Ribet primes.} 
Examples arise from primes $\ell$ so that $\ell \equiv 1 \mod{p}$ and $\ov{\rho}(\Frob_{\ell})=1$. The latter condition is evidently verified for a positive density of primes by Cebotarev's density theorem.

Our construction relies on the following key. 
\begin{prop}\label{prop:lr_ord}
Assume that $p\geq 7$ and $\bar{\rho} \vert _{G_{\BQ(\zeta_p)}}$ is irreducible. 
Let $T$ be a finite set of Ribet primes. 
Then there exists an elliptic newform $h$, giving rise to $\bar{\rho}$, such that the $p$--adic representation $\rho$ associated to $h$ is ordinary at $p$, ramified for $\ell\in T$, and unramified for $\ell' \nmid N \prod\limits_{\ell\in T}\ell$.
\end{prop}

\begin{proof}
 By \cite[Thm.~3]{DT}, there exists a newform $g$, of level dividing $N \prod\limits_{\ell\in T}\ell$, yielding a representation $\rho$ reducing to $\bar{\rho}$ and such that $\rho|_{G_\ell}$ is ramified for $\ell \in T$.
 
As for $p$-ordinarity, we resort to \cite[Thm. 6.1.9 and Lem. 6.1.6]{BLGG}. Note that Lemma 6.1.6 applies to a residual representation of the form 

\begin{equation}
\bar{\rho}|_{G_{\BQ_p}} \simeq \begin{pmatrix}
\bar{\psi_1} &*\\
&\bar{\psi_1} \omega^{-1}
\end{pmatrix},\ 
\end{equation}
rather than the expression \eqref{eq:res}. We apply it to $\bar{\rho}_1 = \omega^{-1} \bar{\rho}$ where $\bar{\rho}$ verifies \eqref{eq:res}. Since $\bar{\rho}_1$ has a Barsotti-Tate lift, so does $\bar{\rho}$.
\end{proof}
\vspace{2mm}

We can now proceed as in section \ref{s:ex}, by choosing an increasing sequence $T_m$ of finite sets of Ribet primes, $T_m=T_{m-1}\cup \{\ell_m\}$, defining ramification sets $S_m$. For all $m$, this defines $R_m^\infty$, the ordinary deformation ring associated to $\rho_m|_G$, $\rho_m$ being the $p$--adic representation associated to the form $h$ of Proposition~\ref{prop:lr_ord}. The arguments of \S\ref{s:ex} extend directly. We deduce the following.

\begin{thm}
Let $p\geq 7$ be a prime. 
Assume $\bar{\rho}$, verifying the conditions in Theorem~\ref{pre}, comes from an elliptic newform $f$ of weight $2$. 

With $R_m^\infty$ as above, assume $R_m^\infty$ is Noetherian for all~$m$. 
Then
 $$R_m^\infty \cong \BZ_p[\![X_1,\ldots X_{s_m}]\!], \text{ $s_m\ge m+1$}.$$
\end{thm}
\vskip2mm
\subsection{The case $\kappa=p+1$} 
This subsection is similar to \S\ref{ss:wt-2}.

Assume that $p\ge 5$. 
We start with a newform $f$ on $\Gamma_1(N)$. If $\ell \nmid Np$ is a prime so that 
$\ell\equiv 1\mod{p}$ and $\ov{\rho}(\Frob_\ell)=1$, Diamond and Taylor \cite[Thm~C]{DT2} show the existence of a newform $h$, giving rise to $\bar{\rho}$, in $S_{p+1}(\Gamma_1(N\ell))$, and new at $\ell$.

\vspace{2mm}

\begin{prop}
If $\bar{\rho}$ is ordinary, then $h$ as above can be chosen to be ordinary at $p$. 
\end{prop}
\begin{proof}
This follows from the result of Diamond-Taylor \cite{DT2} and the results in \cite{BLGG}, exactly as in the proof of Proposition~\ref{prop:lr_ord}.
\end{proof}

Now, choosing a suitable sequence $\ell_1,\ldots, \ell_m,\ldots$ of Ribet primes, we obtain inductively forms $h_m$, of level dividing $N\ell_1\cdots \ell_m$ and new at $\ell_m$. Again, we obtain ramification sets $S_m$ with $S_m\subset S_{m-1} \cup\{\ell_m\}$ and $\ell_m\in S_m$. The same argument, with $R_m^\infty$ defined by the ramification set~$S_m$, yields in this case:

\begin{thm}
Let $p\geq 7$ be a prime. 
Assume $\bar{\rho}$, verifying the conditions in Theorem~\ref{pre}, comes from an elliptic newform of weight $p+1$.

With $R_m^\infty$ as above, assume $R_m^\infty$ is Noetherian for all~$m$. Then
 $$R_m^\infty \cong \BZ_p[\![X_1,\ldots X_{s_m}]\!], \text{ $s_m\ge m+1$}.$$

\end{thm}

\vspace{2mm}

\begin{remark}
\noindent
\begin{itemize}
\item[(1)] The arguments in this section extend to a totally real field unramified at $p$.  
\item[(2)] It is likely that the results extend to representations $\bar{\rho}$ \textit{ordinary of weight} $k\ge 2$ (obvious definition). We leave this to the interested reader.
\end{itemize}
\end{remark}

\end{document}